\documentclass[11pt,a4paper]{article}

\usepackage[utf8]{inputenc}
\usepackage{amsmath,amsfonts,amssymb,amsthm}
\usepackage{geometry}
\usepackage{tikz}
\usetikzlibrary{positioning, shapes.geometric, calc, shadows} 
\usepackage{hyperref}
\usepackage{xcolor}
\usetikzlibrary{shadows, positioning, shapes.geometric}
\newtheorem{theorem}{Theorem}
\newtheorem{lemma}{Lemma}
\newtheorem{definition}{Definition}
\newtheorem{remark}{Remark}
\title{\textbf{A Structural Lower Bound for the Ramsey Number of the Petersen Graph}}
\author{Yuan Si \\
	\texttt{yuan\_si@aliyun.com}}
\date{}

\begin{document}
	
	\maketitle
	
	\begin{abstract}
		In this note, we investigate the diagonal Ramsey number of the Petersen graph, denoted by \(R(\mathcal{P}, \mathcal{P})\). To the best of our knowledge, the dynamic survey \textit{Small Ramsey Numbers} (Revision \#18, 2026) records the lower bound as 19, established by Harborth and Krause. We construct a novel 2-coloring on the complete graph \(K_{20}\) that avoids any monochromatic Petersen graph. By defining the blue graph as the complement of the red graph, we utilize a specific join graph structure \(K_{9,9} \vee \overline{K_2}\). We rigorously prove its Petersen-free property, thereby improving the previously documented lower bound to \(R(\mathcal{P}, \mathcal{P}) \ge 21\). Furthermore, we demonstrate the limit of this structural family by showing that the natural extension \(K_{9,9} \vee \overline{K_3}\) contains a Petersen subgraph.
	\end{abstract}
	
	\section{Introduction}
	Ramsey theory studies the conditions under which order must appear. The diagonal Ramsey number \(R(G, G)\) for a graph \(G\) is the minimum integer \(N\) such that any 2-coloring of the edges of \(K_N\) contains a monochromatic copy of \(G\). 
	
	The Petersen graph, denoted by \(\mathcal{P}\), is one of the most celebrated graphs in graph theory. It is a 3-regular graph with 10 vertices, 15 edges, girth 5 (i.e., it contains no triangles and no 4-cycles), and it is non-bipartite. The exact value of \(R(\mathcal{P}, \mathcal{P})\) remains an open problem. The authoritative dynamic survey \textit{Small Ramsey Numbers} by Radziszowski \cite{Radziszowski2026} records in Section 5.13, page 47, that \(R(\mathcal{P}, \mathcal{P}) \ge 19\), a bound established by Harborth and Krause \cite{Harborth2006} using distance colorings.
	
	In this note, we propose a purely structural construction based on complete bipartite graphs, graph complements, and join operations. Our construction improves the documented lower bound to 21. We provide rigorous combinatorial proofs and analyze the limit of this approach.
	
	\section{Preliminaries}
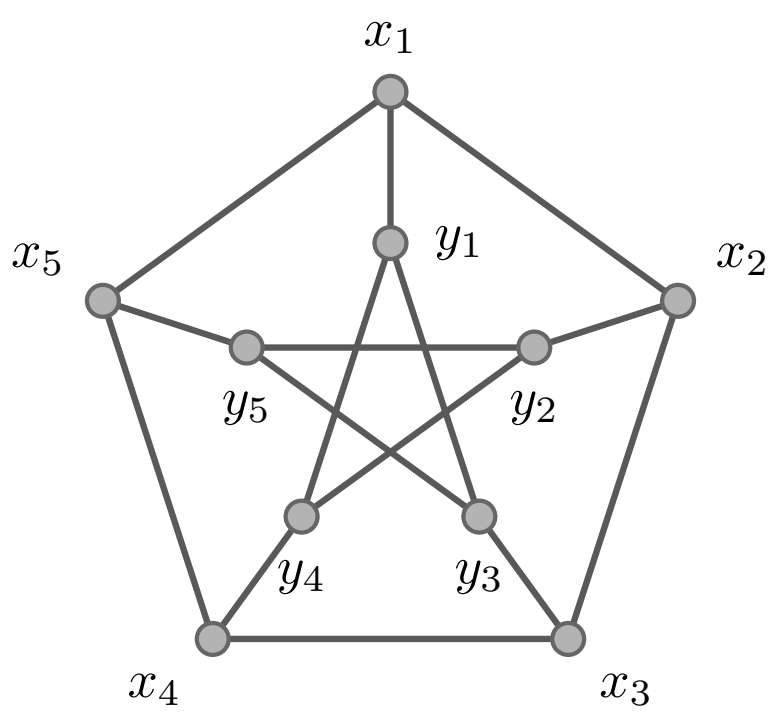
\begin{figure}[htbp]
	\centering
\begin{tikzpicture}[scale=1]
	% 统一定义顶点和边的样式
	\tikzset{
		v/.style={circle, draw=gray!80!black, fill=gray!60, minimum size=6pt, inner sep=0pt, thick},
		e/.style={line width=1.2pt, draw=gray!70!black}
	}
	
	% 1. 计算外圈 5 个顶点的坐标
	\foreach \a in {90, 162, 234, 306, 18} {
		\coordinate (O\a) at (\a:2);
	}
	% 计算内圈 5 个顶点的坐标
	\foreach \a in {90, 162, 234, 306, 18} {
		\coordinate (I\a) at (\a:1);
	}
	
	% 2. 绘制边（稍微加粗）
	% 外圈 5 边形
	\draw[e] (O90) -- (O162) -- (O234) -- (O306) -- (O18) -- (O90);
	% 内圈 5 角星
	\draw[e] (I90) -- (I234) -- (I18) -- (I162) -- (I306) -- (I90);
	% 5 条辐条
	\draw[e] (O90) -- (I90);
	\draw[e] (O162) -- (I162);
	\draw[e] (O234) -- (I234);
	\draw[e] (O306) -- (I306);
	\draw[e] (O18) -- (I18);
	
	% 3. 覆盖绘制实心灰色顶点，按顺时针排序并定制内圈标签方向
	% 循环参数分别为：节点绝对角度/序号/内圈标签的局部角度
	\foreach \ang/\idx/\labang in {90/1/0, 18/2/-90, 306/3/-90, 234/4/-90, 162/5/-90} {
		
		% 外圈节点：标签向外放置（沿 \ang 方向）
		\node[v, label={[text=black, font=\small]\ang:\(x_\idx\)}] at (O\ang) {};
		
		% 内圈节点：标签按定制的 \labang 方向放置，并拉近距离
		\node[v, label={[text=black, font=\small, label distance=1pt]\labang:\(y_\idx\)}] at (I\ang) {};
	}
\end{tikzpicture}
	\caption{The Petersen graph \(\mathcal{P}\).}
	\label{fig:petersen}
\end{figure}

	Let \(\mathcal{P}\) denote the Petersen graph. A standard drawing of \(\mathcal{P}\) is shown in Figure \ref{fig:petersen}, consisting of an outer 5-cycle, an inner 5-cycle, and 5 spokes connecting the corresponding vertices.
	
Throughout this note, we employ several standard graph operations. For two graphs \(G\) and \(H\), their \textbf{union}, denoted by \(G \cup H\), is the graph with vertex set \(V(G) \cup V(H)\) and edge set \(E(G) \cup E(H)\). The \textbf{complement} of a graph \(G\), denoted by \(\overline{G}\), is a graph on the same vertex set \(V(G)\) such that two distinct vertices are adjacent in \(\overline{G}\) if and only if they are not adjacent in \(G\). Furthermore, we utilize the \textbf{join operation} \(\vee\). For disjoint graphs \(G\) and \(H\), \(G \vee H\) is formed by taking the union \(G \cup H\) and adding all edges between the vertices of \(G\) and \(H\).

	In a 2-coloring of \(K_N\), the red graph \(G_R\) and the blue graph \(G_B\) are exact complements of each other; i.e., \(G_B = \overline{G_R}\). Here, \(\overline{K_2}\) denotes the complement of \(K_2\), which is a graph consisting of two isolated vertices.
	
	\section{The Petersen-Free Construction}
	
	\begin{definition}
		Let \(K_{9,9} \vee \overline{K_2}\) be the graph formed by taking a complete bipartite graph with bipartition \(U = \{u_1,u_2, \dots, u_9\}\) and \(W = \{w_1, w_2, \dots, w_9\}\), and joining it with two isolated vertices \(v_1, v_2\) (which form \(\overline{K_2}\)). That is, \(v_1\) and \(v_2\) are not adjacent to each other, and both are connected to every vertex in \(U \cup W\).
	\end{definition}
	
	\begin{lemma}
		\label{lem:no_petersen}
		The graph \(K_{9,9} \vee \overline{K_2}\) does not contain the Petersen graph \(\mathcal{P}\) as a subgraph.
	\end{lemma}
	
	\begin{proof}
		Suppose for the sake of contradiction that \(K_{9,9} \vee \overline{K_2}\) contains a copy of \(\mathcal{P}\). Let \(V(\mathcal{P}) = V_{\mathcal{P}} \subset U \cup W \cup \{v_1, v_2\}\) with \(|V_{\mathcal{P}}| = 10\). 
		Since \(\mathcal{P}\) is non-bipartite and contains 5-cycles, while any subgraph induced by vertices in \(U \cup W\) is a subgraph of the bipartite graph \(K_{9,9}\) and cannot contain odd cycles, \(\mathcal{P}\) cannot be entirely contained in \(U \cup W\).
		
		We consider three cases based on the inclusion of \(v_1, v_2\):
		
		\begin{enumerate}
			\item \textbf{Case 1: \(v_1, v_2 \in V_{\mathcal{P}}\).} 
			Since \(v_1\) and \(v_2\) are not adjacent in \(K_{9,9} \vee \overline{K_2}\), they are not adjacent in \(\mathcal{P}\). Let \(\mathcal{P} - \{v_1, v_2\}\) denote the induced subgraph of \(\mathcal{P}\) obtained by deleting vertices \(v_1, v_2\) and all their incident edges. Since \(v_1, v_2\) are the only vertices outside of \(U \cup W\), the remaining 8 vertices \(V(\mathcal{P} - \{v_1, v_2\})\) must be entirely contained in the bipartite graph \(K_{9,9}\), which means \(\mathcal{P} - \{v_1, v_2\}\) must be bipartite and contain no 5-cycles. We will show this is impossible by proving that \(\mathcal{P} - \{v_1, v_2\}\) always contains a 5-cycle.
			
			Referring to the standard drawing of the Petersen graph in Figure \ref{fig:petersen}, let the outer 5-cycle be  \(\{x_1, x_2, x_3, x_4, x_5\}\) and the inner 5-cycle be \(\{y_1, y_2, y_3, y_4, y_5\}\), where \(x_i\) is connected to \(y_i\) by a spoke.
			
			By the symmetry of \(\mathcal{P}\), we can analyze the positions of \(v_1, v_2\):
			\begin{itemize}
				\item \textbf{Subcase 1a: \(v_1, v_2\) are both in the outer cycle  \(\{x_1, x_2, x_3, x_4, x_5\}\)(or both in the inner cycle \(\{y_1, y_2, y_3, y_4, y_5\}\))} .
				Then the remaining 5 vertices of the other cycle (inner cycle or outer cycle) are completely preserved in \(\mathcal{P} - \{v_1, v_2\}\), directly forming a 5-cycle. Contradiction.
				
				\item \textbf{Subcase 1b: One vertex is in the outer cycle and one is in the inner cycle.}
				Without loss of generality, fix \(v_1 = x_1\). By the symmetry of \(\mathcal{P}\), the remaining four vertices in the inner cycle fall into two symmetry classes relative to \(x_1\): the class \(\{y_2, y_5\}\) and the class \(\{y_3, y_4\}\). Thus, without loss of generality, there are two cases for the position of \(v_2\). If \(v_2 = y_2\) (or symmetrically \(y_5\)), we delete \(x_1\) and \(y_2\). Then \(\mathcal{P} - \{x_1, y_2\}\) still contains the 5-cycle:
					\[ x_5 - y_5 - y_3 - x_3 - x_4 - x_5. \] If \(v_2 = y_3\) (or symmetrically \(y_4\)), we delete \(x_1\) and \(y_3\). Then \(\mathcal{P} - \{x_1, y_3\}\) still contains the 5-cycle:
					\[ y_2 - y_5 - x_5 - x_4 - y_4 - y_2. \]
			
				In both cases, a 5-cycle remains intact.
			\end{itemize}
			Therefore, \(\mathcal{P} - \{v_1, v_2\}\) always contains a 5-cycle, which contradicts the fact that its vertices lie entirely in the bipartite graph \(K_{9,9}\).
			
		 \item \textbf{Case 2: Exactly one of \(v_1, v_2\) is in \(V_{\mathcal{P}}\).} 
		Without loss of generality, assume \(v_1 \in V_{\mathcal{P}}\) and \(v_2 \notin V_{\mathcal{P}}\). Then \(\mathcal{P} - \{v_1\}\) has 9 vertices entirely contained in \(U \cup W\). 
		In the Petersen graph \(\mathcal{P}\), there are two disjoint 5-cycles: the outer cycle \(\{x_1, x_2, x_3, x_4, x_5\}\) and the inner cycle \(\{y_1, y_2, y_3, y_4, y_5\}\). Deleting any single vertex \(v_1\) can destroy at most one of these two 5-cycles. Thus, \(\mathcal{P} - \{v_1\}\) must still contain at least one 5-cycle. 
		However, its vertices lie entirely in \(U \cup W\), which induces the bipartite graph \(K_{9,9}\) and cannot contain any odd cycles (such as a 5-cycle). Contradiction.
			
			\item \textbf{Case 3: Neither \(v_1\) nor \(v_2\) is in \(V_{\mathcal{P}}\).} 
			Then \(V_{\mathcal{P}} \subset U \cup W\). This is the case already excluded at the beginning of the proof.
		\end{enumerate}
		Thus, no such \(\mathcal{P}\) exists in \(K_{9,9} \vee \overline{K_2}\).
	\end{proof}
	
\begin{figure}[htbp]
		\centering  % 必须在 figure 内，tikzpicture 外
\begin{tikzpicture}[
	scale=1.2,
	% 定义红色连线样式
	red_edge/.style={line width=1.5pt, draw=red!80!black},
	% 定义蓝色基础样式（仅设定颜色和双线，不设定粗细）
	blue_base/.style={draw=blue!80!black, double=blue!20!white},
	% 定义方形集合样式（带阴影） - 尺寸稍微增大以容纳标签
	set_blob/.style={
		rectangle, thick, draw=red!80!black, fill=red!15, 
		minimum width=3cm, minimum height=3cm, rounded corners=4pt,
		drop shadow={shadow xshift=2pt, shadow yshift=-2pt, opacity=0.35}
	},
	% 定义灰色实心顶点样式（带阴影）
	vertex/.style={
		circle, draw=gray!80!black, fill=gray!60, 
		minimum size=10pt, inner sep=0pt, thick,
		circular drop shadow={shadow xshift=1pt, shadow yshift=-1pt, opacity=0.4}
	}
	]
	
	% 1. 先定义底层的不可见坐标点
	\coordinate (U_pos) at (-3, 0);
	\coordinate (W_pos) at (3, 0);
	\coordinate (v1_pos) at (-1.5, 2.8);
	\coordinate (v2_pos) at (1.5, 2.8);
	
	% 提前计算框内 9 个点的坐标（绕着中心排成圆，半径缩小为 0.85cm，确保标签不超出边框）
	\foreach \i in {1,...,9} {
		\pgfmathsetmacro{\ang}{90-40*(\i-1)}
		\path (U_pos) + (\ang:0.85cm) coordinate (U_\i);
		\path (W_pos) + (\ang:0.85cm) coordinate (W_\i);
	}
	
	% 2. 绘制底层的蓝色连线（从中心点出发，会被后面的框和点遮挡）
	% 【v1 连向 U 的线】
	\draw[blue_base, line width=1.5pt, double distance=6pt, bend left=15] (v1_pos) to (U_pos);
	% 【v1 连向 W 的线】
	\draw[blue_base, line width=1.5pt, double distance=6pt, bend right=10] (v1_pos) to (W_pos);
	% 【v2 连向 W 的线】
	\draw[blue_base, line width=1.5pt, double distance=6pt, bend right=15] (v2_pos) to (W_pos);
	% 【v2 连向 U 的线】
	\draw[blue_base, line width=1.5pt, double distance=6pt, bend left=10] (v2_pos) to (U_pos);
	
	% 3. 绘制底层的红色 K_2 边
	\draw[red_edge] (v1_pos) -- (v2_pos);
	
	% 4. 绘制 U 和 W 之间的粗蓝双线及 K_{9,9}
	\draw[blue_base, line width=2.5pt, double distance=25pt] (U_pos) -- (W_pos) 
	node[midway, text=blue!80!black, font=\bfseries\large] {\(K_{9,9}\)};
	
	% 5. 绘制顶层：覆盖红色方块（带下方标签）
	\node[set_blob, label={[text=red!80!black, font=\bfseries\Large]270:\(K_9\)}] (U) at (U_pos) {};
	\node[set_blob, label={[text=red!80!black, font=\bfseries\Large]270:\(K_9\)}] at (W_pos) {};
	
	% 6. 在框内绘制 K_9 的边和顶点
	% 为了避免密集的线糊成一团，框内部的红色连线画细一点
	\foreach \i in {1,...,9} {
		\foreach \j in {1,...,9} {
			\ifnum\i<\j
			\draw[red_edge, line width=0.6pt] (U_\i) -- (U_\j);
			\draw[red_edge, line width=0.6pt] (W_\i) -- (W_\j);
			\fi
		}
	}
	% 画小圆点，放在红边上面，并使用独立阴影，同时加上黑色的小号顶点标签
	% 微调：label distance=1pt, inner sep=0pt 让标签紧紧贴靠顶点
	\foreach \i in {1,...,9} {
		\pgfmathsetmacro{\ang}{90-40*(\i-1)}
		\node[vertex, minimum size=4pt, circular drop shadow={shadow xshift=0.5pt, shadow yshift=-0.5pt, opacity=0.4}, label={[text=black, font=\scriptsize, label distance=1pt, inner sep=0pt]\ang:\(u_\i\)}] at (U_\i) {};
		\node[vertex, minimum size=4pt, circular drop shadow={shadow xshift=0.5pt, shadow yshift=-0.5pt, opacity=0.4}, label={[text=black, font=\scriptsize, label distance=1pt, inner sep=0pt]\ang:\(w_\i\)}] at (W_\i) {};
	}
	
	% 7. 绘制黑色实心顶点 v1, v2 及 K_2 文本
	\node[vertex, label={[text=black, font=\bfseries]90:\(v_1\)}] at (v1_pos) {};
	\node[vertex, label={[text=black, font=\bfseries]90:\(v_2\)}] at (v2_pos) {};
	\node[text=red!80!black, font=\bfseries\large] at (0, 3.2) {\(K_2\)};
	
	% 8. 图例
	\node[draw=none, rectangle, anchor=west, text=black, font=\small, align=left, rounded corners, fill=gray!10, inner sep=8pt] at (-4.5, -2.5) {
		\tikz{\draw[red_edge] (0,0) -- (0.6,0);} Red graph \(G_R = K_9 \cup K_9 \cup K_2\) \\[1ex]
		\tikz{\draw[blue_base, line width=1.5pt, double distance=4pt] (0,0) -- (0.6,0);} Blue graph \(G_B = K_{9,9} \vee \overline{K_2}\)
	};
	
\end{tikzpicture}
\caption{The 2-coloring of \(K_{20}\). The red graph consists of two \(K_9\) components (rectangular blocks) and one \(K_2\) edge. The blue graph (thick double lines) is its exact complement, containing all cross-edges between the sets.}
\label{fig:construction}
\end{figure}
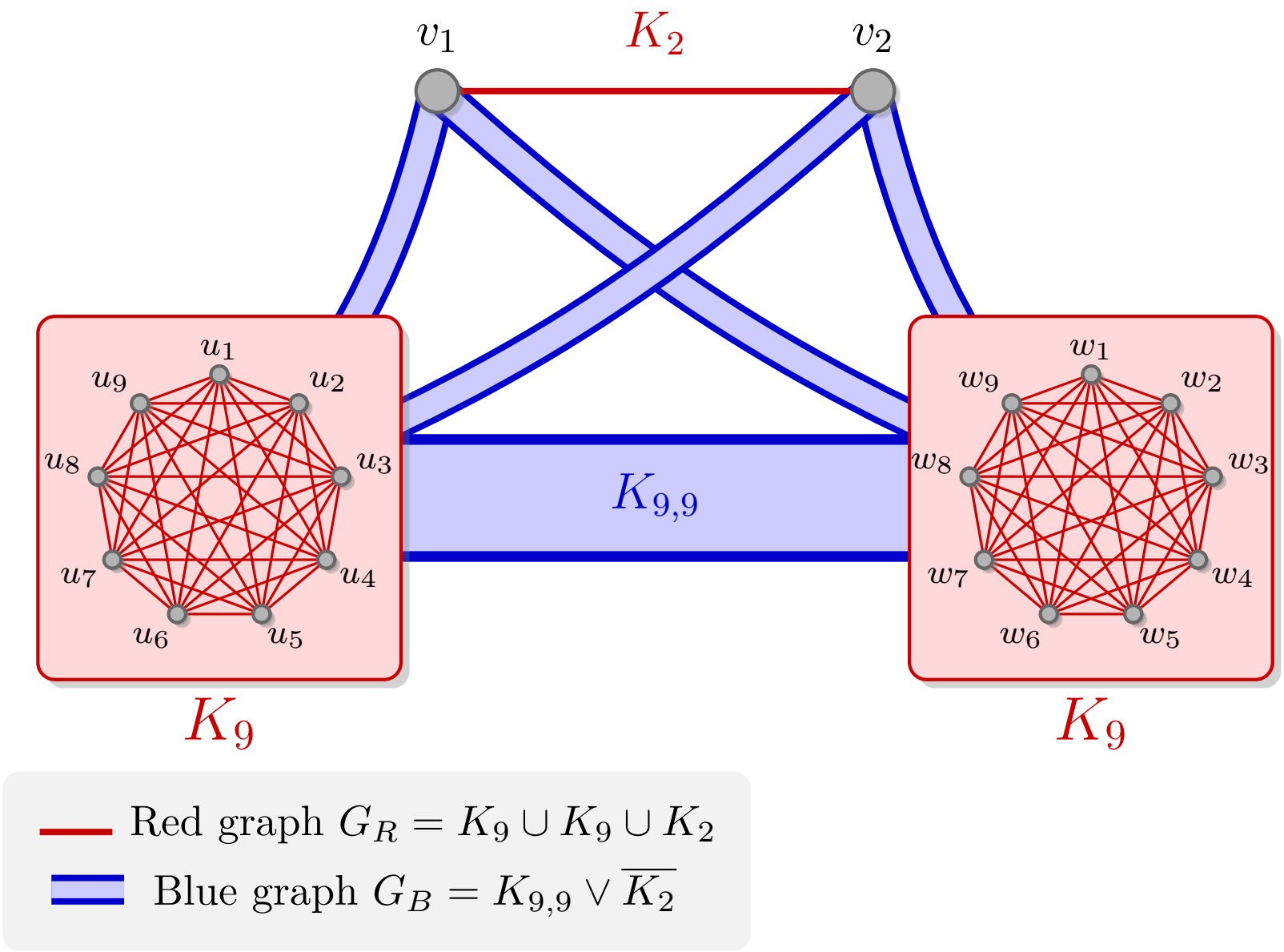

\begin{theorem}
	\label{thm:main}
	\(R(\mathcal{P}, \mathcal{P}) \ge 21\).
\end{theorem}
\begin{proof}
	We construct a 2-coloring of the edges of \(K_{20}\) as follows:
	\begin{itemize}
		\item \textbf{Red graph \(G_R\):} Consists of two disjoint copies of \(K_9\) and one copy of \(K_2\). The total number of vertices is \(9 + 9 + 2 = 20\). The Petersen graph has 10 vertices and is connected; hence a graph whose largest connected component has size 9 cannot contain it.
		\item \textbf{Blue graph \(G_B\):} The blue graph is the complement of the red graph. The complement of \(K_9 \cup K_9 \cup K_2\) is exactly \(K_{9,9} \vee \overline{K_2}\). By Lemma \ref{lem:no_petersen}, \(G_B\) contains no Petersen graph.
	\end{itemize}
	Since neither the red nor the blue subgraph contains a monochromatic Petersen graph, the Ramsey number is at least 21.
\end{proof}

	\section{The Limit of the Construction}
	
	It is natural to ask whether this structural family can be extended to \(K_{21}\). Suppose we set the red graph as \(K_9 \cup K_9 \cup K_3\). The complement of the red graph becomes \(K_{9,9} \vee \overline{K_3}\). 
	
	\begin{theorem}
		\label{thm:limit}
		The graph \(K_{9,9} \vee \overline{K_3}\) contains a Petersen graph.
	\end{theorem}
	
	\begin{proof}
	Let the vertices of \(\overline{K_3}\) be \(v_1, v_2, v_3\), and recall that \(K_{9,9}\) has the bipartition \(U = \{u_1, u_2, \dots, u_9\}\) and \(W = \{w_1, w_2, \dots, w_9\}\).  We explicitly construct a subgraph using the following 10 vertices:
		\[
		V_{\mathcal{P}} = \{v_1, v_2, v_3, u_1, u_2, u_3, w_1, w_2, w_3, w_4\}.
		\]
		The 15 edges are given as follows:
		\begin{itemize}
			\item \textbf{Edges incident to \(\overline{K_3}\):} \(v_1u_1, v_1w_1, v_1w_2;\quad v_2u_2, v_2w_1, v_2w_3;\quad v_3u_3, v_3w_2, v_3w_3.\)
			\item \textbf{Edges between \(U\) and \(W\):} \(u_1w_3, u_1w_4;\quad u_2w_2, u_2w_4;\quad u_3w_1, u_3w_4.\)
		\end{itemize}
		
		To verify this is isomorphic to \(\mathcal{P}\), observe that these edges decompose into two 5-cycles plus a perfect matching:
		\begin{itemize}
			\item \textbf{Outer 5-cycle :} \(v_1 - w_1 - v_2 - w_3 - u_1 - v_1\).
			\item \textbf{Inner 5-cycle :} \(v_3 - u_3 - w_4 - u_2 - w_2 - v_3\).
			\item \textbf{Perfect Matching:} \(v_1w_2,\ w_1u_3,\ v_2u_2,\ w_3v_3,\ u_1w_4\).
		\end{itemize}
	\end{proof}
	
\begin{figure}[htbp]
	\centering
	\begin{tikzpicture}[scale=1]
		% 统一定义顶点和边的样式
		\tikzset{
			v/.style={circle, draw=gray!80!black, fill=gray!60, minimum size=6pt, inner sep=0pt, thick},
			e/.style={line width=1.2pt, draw=gray!70!black}
		}
		
		% 1. 计算外圈 5 个顶点的坐标
		\foreach \a in {90, 162, 234, 306, 18} {
			\coordinate (O\a) at (\a:2);
		}
		% 计算内圈 5 个顶点的坐标
		\foreach \a in {90, 162, 234, 306, 18} {
			\coordinate (I\a) at (\a:1);
		}
		
		% 2. 绘制边（稍微加粗）
		% 外圈 5 边形
		\draw[e] (O90) -- (O162) -- (O234) -- (O306) -- (O18) -- (O90);
		% 内圈 5 角星
		\draw[e] (I90) -- (I234) -- (I18) -- (I162) -- (I306) -- (I90);
		% 5 条辐条 (Perfect Matching)
		\draw[e] (O90) -- (I90);
		\draw[e] (O162) -- (I162);
		\draw[e] (O234) -- (I234);
		\draw[e] (O306) -- (I306);
		\draw[e] (O18) -- (I18);
		
		% 3. 覆盖绘制实心灰色顶点，按顺时针排序并定制内圈标签方向
		% 严格按照你提供的边集映射：
		% 外圈: v1 - w1 - v2 - w3 - u1
		% 内圈: w2 - u3 - u2 - v3 - w4
		\foreach \ang/\xlabel/\ylabel/\labang in {
			90/\(v_1\)/\(w_2\)/0, 
			162/\(w_1\)/\(u_3\)/-90, 
			234/\(v_2\)/\(u_2\)/-90, 
			306/\(w_3\)/\(v_3\)/-90, 
			18/\(u_1\)/\(w_4\)/-90} {
			
			% 外圈节点：标签向外放置
			\node[v, label={[text=black, font=\small]\ang:\xlabel}] at (O\ang) {};
			
			% 内圈节点：标签按定制的 \labang 方向放置，完美避开连线
			\node[v, label={[text=black, font=\small, label distance=1pt]\labang:\ylabel}] at (I\ang) {};
		}
	\end{tikzpicture}
	\caption{Isomorphic mapping of the constructed Petersen subgraph in \(K_{9,9} \vee \overline{K_3}\).}
	\label{fig:petersen_isomorphic}
\end{figure}
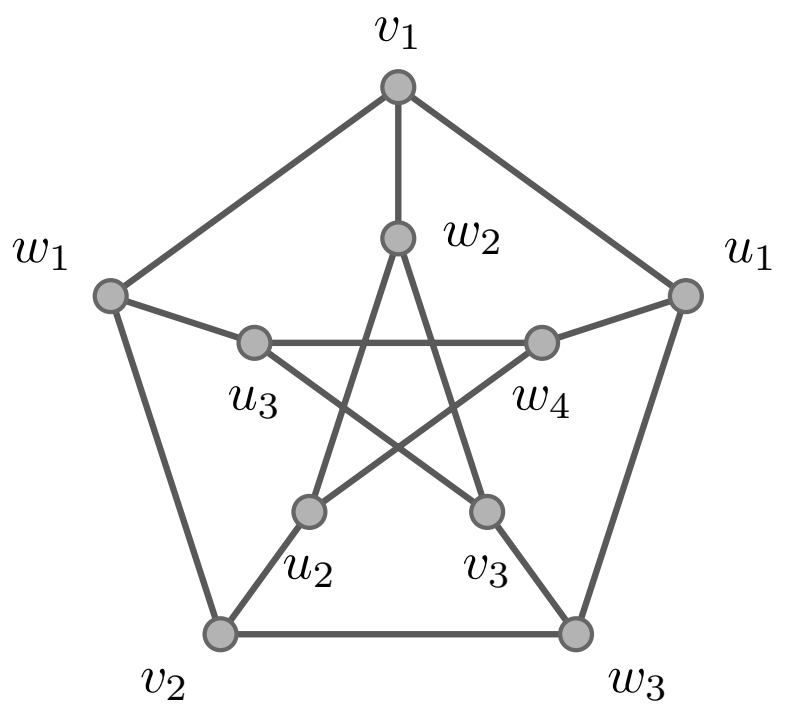
	
	\begin{remark}
		Theorem \ref{thm:limit} demonstrates that the specific complementary construction \(K_{9,9} \vee \overline{K_3}\) contains a Petersen subgraph. It does not establish an upper bound for \(R(\mathcal{P}, \mathcal{P})\); the exact value could still be greater than 21.
	\end{remark}

\end{document}